\documentclass{amsart}

\usepackage{times}
\usepackage{amssymb,amsmath,enumerate,amsthm}
\usepackage[all]{xypic}

\newtheorem{theorem}{Theorem}[section]

\newtheorem{axiom}[theorem]{Axiom}

\begin{document} 

\title{Bivariant Versions of Algebraic Cobordism} 

\author{Rui Miguel Saramago} 
\address{Departamento de Matem\'atica, Instituto Superior T\'ecnico, Av. Rovisco Pais, 1049-001 Lisboa, Portugal} 
\email{saramago@math.tecnico.ulisboa.pt} 


\begin{abstract}  
We define four distinct oriented bivariant theories associated with algebraic cobordism in its two versions (the axiomatic $\Omega$ and the geometric $\omega$), when applied to quasi-projective varieties over a field $k$. Specifically, we obtain contravariant analogues of the algebraic bordism group $\Omega_*(X)$ and the double point bordism group $\omega_*(X)$, for $X$ a quasi-projective variety, and covariant analogues of the algebraic cobordism ring $\Omega^*(X)$ and the double point cobordism ring $\omega^*(X)$, for $X$ a smooth variety. When the ground field has characteristic zero, we use the universal properties of algebraic cobordism in order to obtain correspondences between these oriented bivariant theories.
\end{abstract} 

\maketitle
\tableofcontents

\section{Introduction}

This work aims at positioning the algebraic and double point cobordism theories of Levine-Morel \cite{lm} and Levine-Pandharipande \cite{lp} in a bivariant theory setting. This project was directly inspired by the survey by Sch\"urman and Yokura in \cite{sy}, where the problem is presented. The main interest in doing this comes from the added richness provided by an (oriented) bivariant theory: in particular, any such $\mathbb{B^*} (X \xrightarrow{f} Y)$, defined in a category where the objects have direct geometric meaning, gives rise to a covariant theory $\mathbb{B_*} (X) := \mathbb{B^{-*}} (X \xrightarrow{f} pt)$ and to a contravariant theory $\mathbb{B^*} (X) := \mathbb{B^*} (X \xrightarrow{id} X)$. Providing convenient bivariant theories, one can then in principle obtain \textit{mirror} theories for algebraic cobordism and algebraic bordism. For example, algebraic cobordism (in its first iteration from \cite{lm}) is a universal oriented cohomology theory (defined e.g. on a category of quasi-projective varieties), usually obtained by "Poincar\'e duality" from algebraic bordism, a universal oriented Borel-Moore homology theory; considering it as part of a bivariant theory, which we call $\mathbb{OB}2$, the corresponding covariant part will be a new oriented Borel-Moore homology theory, different from algebraic bordism.

The construction of each of the bivariant theories that will be presented starts from $\mathbb{OB}$, defined by Sch\"urman and Yokura in \cite{sy}, and from one analogue, called $\mathbb{OB'}$, introduced in section 6. $\mathbb{OB}$ was in turn inspired by the group of cobordism cycles $\mathcal{Z}_*$ of Levine and Morel from \cite{lm}. Some of our bivariant theories will be obtained by modding $\mathbb{OB}$ and $\mathbb{OB'}$ by relations (of a geometric nature), in the same way as algebraic bordism is constructed from $\mathcal{Z}_*$ by modding out in succession the dimension, section and formal group law axioms.

This work starts by recalling, in section 2, the definition and properties of Fulton-MacPherson's bivariant theories that will be used later on. Section 3 presents algebraic bordism and cobordism, when applied resp. to quasi-projective and to smooth varieties over a field $k$. If the characteristic of $k$ is $0$, or if resolution of singularities is at hand, these have precise universal properties among homology and cohomology theories of the same kind. 

In section 4, we recall the definition of the universal oriented bivariant theory $\mathbb{OB}$ for quasi-projective varieties, and in section 5 we use this theory to define $\mathbb{OB}1$, a bivariant theory whose covariant part is algebraic bordism $\Omega_*$ for quasi-projective varieties. This definition is done through a generalization of the concept of fundamental class for a smooth variety, that is, a class $[X \xrightarrow{f} Y] \in \mathbb{OB}(X \xrightarrow{f} Y)$ that, for $Y$ a point, is the same as the known definition of $[X]$. We describe the contravariant theory $\mathbb{OB}1(X \xrightarrow{id} X)$ for smooth varieties. In section 6, a new bivariant theory is introduced, based on a "dual" version $\mathbb{OB'}$ of $\mathbb{OB}$, which we call $\mathbb{OB}2$, and whose contravariant part is algebraic cobordism $\Omega^*$ for smooth varieties. We proceed by deducing the covariant part of $\mathbb{OB}2$ for quasi-projective varieties.

In section 7 we recall Levine-Pandharipande's definition of double point bordism, and the associated double-point cobordism, and in sections 8 and 9 we further the program from sections 5 and 6. Namely, we define a bivariant theory $\mathbb{OB}3$ whose covariant part $\mathbb{OB}3 (X \xrightarrow{c} pt)$ is double-point bordism, and present the corresponding contravariant theory for smooth varieties; and a bivariant theory  $\mathbb{OB}4$ whose contravariant part, for smooth varieties, is double-point cobordism, showing the corresponding covariant theory for quasi-projective varieties.

The constructions presented are formal and depend deeply on the definition of bivariant theory. The last, brief section tries to connect the four previous bivariant theories, using the universal properties of the theories of Levine-Morel and Levine-Pandharipande.

A remark on notation: the original definition by Levine-Morel denotes algebraic cobordism by $\Omega$, which encompasses both $\Omega_*$ and $\Omega^*$, and the same happens with double point cobordism $\omega$. We decided to call $\Omega_*$ (the quotient of the group of cobordism cycles by the axioms of dimension, section and formal group law) \textit{algebraic bordism}, and save the term \textit{algebraic cobordism} for $\Omega^*$, the "dual" of $\Omega_*$. We follow a similar approach for $\omega$: $\omega_*$ is for us \textit{double point bordism} and $\omega^*$ is \textit{double point cobordism}.

\section{Oriented Bivariant Theories}

This section follows Fulton and MacPerson's terminology from \cite{fm}. \\

Let $\mathcal{C}$ be a category with a final object and with all fiber products. Distinguish a class of morphisms, which we call {\it confined maps}, and a class of fiber squares, called 
{\it confined squares} (which are sometimes also called {\it independent squares}). The class of confined maps must contain all identity maps and be closed under composition and base change. The class of confined squares must satisfy two conditions:

\begin{itemize}
  \item Any square of one the following forms is a confined square: \\ \\
  \xymatrix{
  X \ar[r]^{id_X} \ar[d]_{f} & X \ar[d]^{f} & & & X \ar[r]^{f} \ar[d]_{id_X} & Y \ar[d]^{id_Y}\\
  Y \ar[r]^{id_Y}  & Y & & & X \ar[r]^{f}  & Y } \\ \\
  \item If the two inside squares of any diagram of one of the following forms are confined squares, than so is the respective outside square: \\ \\
  \xymatrix{
  X \ar[r]^{f} \ar[d]_{h} & X' \ar[r]^{f'} \ar[d]_{h'}  & X'' \ar[d]_{h''} && & X \ar[r]^{f} \ar[d]_{a} & X' \ar[d]^{a'}\\
  Y \ar[r]_{g}  & Y' \ar[r]_{g'} & Y'' & & & Y \ar[r]^{g}  \ar[d]_{b} & Y' \ar[d]^{b'} \\
  & & & & & Z \ar[r]^{h}  & Z'}
\end{itemize}

\vspace{0.5cm}
The example that we will use throughout this work is that of quasi-projective varieties over a field of characteristic $0$, with proper morphisms as the confined maps, and all fiber squares as the confined squares. Since we will look at ways to produce bivariant analogues of algebraic cobordism or bordism, we will need resolution of singularities; we could choose a category with different objects, e.g. by changing the characteristic of the base field, provided we require that resolutions of singularities exist for any of the varieties in the category.

We will define graded bivariant theories on a category $\mathcal{C}$ as above. A graded bivariant theory $\mathbb{B^*}$ on $\mathcal{C}$ assigns to each morphism $X \xrightarrow{f} Y$ in the category a graded abelian group, denoted by $\mathbb{B^*} (X \xrightarrow{f} Y)$. On these graded abelian groups we define three operations:
\begin{itemize}
  \item A $\mathbb{Z}$- bilinear {\bf product} $\bullet : {\mathbb B}^i (X \xrightarrow{f} Y) \otimes {\mathbb B}^j (Y \xrightarrow{g} Z) \rightarrow {\mathbb B}^{i+j} (X \xrightarrow{gf} Z)$, \\
  defined for any $X \xrightarrow{f} Y \xrightarrow{g} Z$. \\
 \item A $\mathbb{Z}$- linear {\bf pushforward} $f_* : {\mathbb B}^i (X \xrightarrow{gf} Z) \rightarrow {\mathbb B}^{i} (Y \xrightarrow{g} Z)$,\\
  defined for any  $X \xrightarrow{f} Y \xrightarrow{g} Z$ with $f$ confined. \\
\item A $\mathbb{Z}$- linear {\bf pullback} $g^* : {\mathbb B}^i (X' \xrightarrow{h'} Y') \rightarrow {\mathbb B}^{i} (X \xrightarrow{h} Y)$, \\
  defined for any confined square  \xymatrix{
  X \ar[r]^{f} \ar[d]_{h} & X' \ar[d]^{h'} \\
  Y \ar[r]^{g}  & Y'  }. 
\end{itemize}

The operations are defined for each degree (or pair of degrees, for the product). Moreover, they must satisfy the following axioms, which we describe at length:
\begin{itemize}
  \item The product is associative: $(a \bullet b) \bullet c = a \bullet (b \bullet c)$ for any $a \in {\mathbb B}^i (X \xrightarrow{f} Y)$, $b \in {\mathbb B}^j (Y \xrightarrow{g} Z)$ and $c \in {\mathbb B}^k (Z \xrightarrow{h} W)$ (given any $X \xrightarrow{f} Y \xrightarrow{g} Z \xrightarrow{h} W$). \\
  \item The pushforward is functorial: $(id_X)_* = id_{{\mathbb B}^i (X \xrightarrow{g} Z)}$ for any $Y \xrightarrow{g} Z$ \\
  and $(f_2 f_1)_* = {f_2}_* {f_1}_*$ for any $X \xrightarrow{f_1} X' \xrightarrow{f_2} Y \xrightarrow{g} Z$ with $f_1$ and $f_2$ confined. \\
  \item The pullback is functorial: $(id_Y)* = id_{{\mathbb B}^i (X \xrightarrow{h} Y)}$ for any square \xymatrix{
  X \ar[r]^{id_X} \ar[d]_{h} & X \ar[d]^{h} \\
  Y \ar[r]^{id_Y}  & Y } \\
  and $(g_2 g_1)^* = g_1^* \ g_2^* \ $ for any  \xymatrix{
  X \ar[r]^{f_1} \ar[d]_{h} & X' \ar[r]^{f_2} \ar[d]_{h'}  & X'' \ar[d]_{h''} \\
  Y \ar[r]_{g_1}  & Y' \ar[r]_{g_2} & Y''} \\ with both inner squares confined. \\
  \item The product commutes with the pushforward: $(f_* a) \bullet b = f_*(a \bullet b)$ for any $a \in {\mathbb B}^i (X \xrightarrow{hf} Y)$ and $b \in {\mathbb B}^j (Y \xrightarrow{g} Z)$ (given any $X \xrightarrow{f} X' \xrightarrow{h} Y \xrightarrow{g} Z$ with f confined). \\
  \item The product commutes with the pullback: $(id_Y^* a) \bullet b = id_Z^*(a \bullet b)$ for any $a \in {\mathbb B}^i (X' \xrightarrow{h'} Y)$ and $b \in {\mathbb B}^j (Y \xrightarrow{g} Z)$ (given any \xymatrix{
  X \ar[r]^{f} \ar[d]_{h} & X' \ar[d]^{h'} \\
  Y \ar[r]^{id_Y} \ar[d]_{g} & Y \ar[d]^{g} \\
  Z \ar[r]^{id_Z}  & Z } \\
   with the top square confined). \\
  \item The pushforward commutes with the pullback: $g_*(f_2^* a) = f_3^*(g'_* a)$ for any $a \in {\mathbb B}^i (X' \xrightarrow{g'} Y')$ (given any \xymatrix{
  X \ar[r]^{f_1} \ar[d]_{g} & X' \ar[d]^{g'} \\
  Y \ar[r]^{f_2} \ar[d]_{h} & Y' \ar[d]^{h'} \\
  Z \ar[r]^{f_3}  & Z' } \\
   with $g$ and $g'$ confined maps and both inner squares also confined). \\

  \item Projection: $g'_*(g^* a \bullet b) = a \bullet g_*b$ for any $a \in {\mathbb B}^i (X \xrightarrow{f} Y)$ and $b \in {\mathbb B}^j (Y' \xrightarrow{hg} Z)$ (given any \xymatrix{
  X' \ar[r]^{g'} \ar[d]_{f'} & X \ar[d]^{f} & \\
  Y' \ar[r]^{g} & Y \ar[r]^{h} & Z} \\
   where the square is confined and both $g$ and $g'$ are also confined). \\
   Note: $h g f' = h f g'$, so it makes sense to apply the pushforward of $g'$.\\
   
\end{itemize}

Yokura and Sch\"urman introduce also the notion of an {\it oriented bivariant theory} \cite{y}. This is a bivariant theory $\mathbb{B^*}$, defined on the category of quasi-projective varieties over a field $k$ (with no restriction on characteristic) and proper morphisms, with all fiber squares as confined maps and all projective maps as confined maps, but now with an attribution of a {\it first Chern class operator} for each morphism $X \xrightarrow{f} Y$ and line bundle $L$ over $X$, which is a degree-one morphism of abelian groups 
$$\widetilde{c_1} (L) : \mathbb{B^*}(X \xrightarrow{f} Y) \rightarrow \mathbb{B}^{*-1}(X \xrightarrow{f} Y)$$
which is compatible with the bivariant theory structure: 

\begin{itemize}
  \item The Chern class operators for different line bundles over the same  variety commute: $\widetilde{c_1} (L_1) \circ \widetilde{c_1} (L_2) = \widetilde{c_1} (L_2) \circ \widetilde{c_1} (L_1) : \mathbb{B^*}(X \xrightarrow{f} Y) \rightarrow \mathbb{B}^{*-2}(X \xrightarrow{f} Y)$ for any $X \xrightarrow{f} Y$ and any two line bundles $L_1$ and $L_2$ over $X$. \\
  Additionally, isomorphic line bundles induce the same first Chern class operator. \\
  \item The Chern class operator is compatible with the pushforward: For any $X \xrightarrow{f} Y \xrightarrow{g} Z$, with $f$ confined, and any line bundle $L$ over $Y$, the following square is commutative: 
  \xymatrix{
  \mathbb{B}^i(X \xrightarrow{gf} Z) \ar[r]^{f_*} \ar[d]_{\widetilde{c_1} (f^*L)} & \mathbb{B}^i(Y \xrightarrow{g} Z) \ar[d]^{\widetilde{c_1} (L)} \\
  \mathbb{B}^{i-1}(X \xrightarrow{gf} Z) \ar[r]^{f_*}  & \mathbb{B}^{i-1}(Y \xrightarrow{g} Z) } \\
  \item The Chern class operator is compatible with the pullback: for any confined square \xymatrix{
  X \ar[r]^{f} \ar[d]_{h} & X' \ar[d]^{h'} \\
  Y \ar[r]^{g}  & Y'  } 
  and any line bundle $L$ over $X'$, the following square is commutative: 
  \xymatrix{
  \mathbb{B}^i(X' \xrightarrow{h'} Y') \ar[r]^{g^*} \ar[d]_{\widetilde{c_1} (L)} & \mathbb{B}^i(X \xrightarrow{h} Y) \ar[d]^{\widetilde{c_1} (f^*L)} \\
  \mathbb{B}^{i-1}(X' \xrightarrow{h'} Y') \ar[r]^{g^*}  & \mathbb{B}^{i-1}(X \xrightarrow{h} Y) } \\
  \item The Chern class operator is compatible with the product: For any $X \xrightarrow{f} Y \xrightarrow{g} Z$ and any line bundles $L_1$ over $X$ and $L_2$ over $Y$, the following squares are commutative: \\
  \xymatrix{
  \mathbb{B}^i(X \xrightarrow{f} Y) \otimes \mathbb{B}^j(Y \xrightarrow{g} Z) \ar[r]^{\ \ \ \ \ \ \ \ \bullet} \ar[d]_{\widetilde{c_1} (L_1) \otimes id} & \mathbb{B}^{i+j}(X \xrightarrow{gf} Z) \ar[d]^{\widetilde{c_1} (L_1)} \\
  \mathbb{B}^{i-1}(X \xrightarrow{f} Y) \otimes \mathbb{B}^j(Y \xrightarrow{g} Z) \ar[r]^{\ \ \ \ \ \ \ \ \ \bullet}  & \mathbb{B}^{i+j-1}(X \xrightarrow{gf} Z) } \\
  \\ \\
  \xymatrix{
  \mathbb{B}^i(X \xrightarrow{f} Y) \otimes \mathbb{B}^j(Y \xrightarrow{g} Z) \ar[r]^{\ \ \ \ \ \ \ \ \bullet} \ar[d]_{id \otimes \widetilde{c_1} (L_2)} & \mathbb{B}^{i+j}(X \xrightarrow{gf} Z) \ar[d]^{\widetilde{c_1} (f^*L_2)} \\
  \mathbb{B}^{i}(X \xrightarrow{f} Y) \otimes \mathbb{B}^{j-1}(Y \xrightarrow{g} Z) \ar[r]^{\ \ \ \ \ \ \ \ \ \bullet}  & \mathbb{B}^{i+j-1}(X \xrightarrow{gf} Z) }
\end{itemize}

Note that, in the above diagrams, Chern class operators denoted by the same symbol may actually be different: for example, in the first diagram for the product compatibility,  $\widetilde{c_1} (L_1)$ has a different meaning in the two down arrows, although both refer to the same line bundle $L_1$ over $X$.

Any bivariant theory $\mathbb{B^*}$ induces naturally a covariant and a contravariant theory through the corresponding associated functors. Define the covariant graded theory $\mathbb{B}_i (X) := \mathbb{B}^{-i} (X \xrightarrow{c} pt)$, which gives a covariant functor if one restricts to confined morphisms (and thus $f : X \rightarrow Y$ induces a group morphism $f_* : \mathbb{B}_i (X) \rightarrow \mathbb{B}_i (Y)$ by using the pushforward). On the other hand, define the contravariant theory $\mathbb{B}^i (X) := \mathbb{B}^{i} (X \xrightarrow{id} X)$, a contravariant functor ($f : X \rightarrow Y$ induces a group homomorphism $f^* : \mathbb{B}^i (Y) \rightarrow \mathbb{B}^i (X)$ by using the pullback associated with the square \xymatrix{
  X \ar[r]^{f} \ar[d]_{id} & Y \ar[d]^{id} \\
  X \ar[r]^{f}  & Y  }, which is confined for all $f : X \rightarrow Y$).
  For each $X$, $\mathbb{B}^i (X)$ has actually a structure of graded ring: the product $\bullet$ of the bivariant theory induces a product $\bullet : \mathbb{B}^{i} (X \xrightarrow{id} X) \otimes \mathbb{B}^{i} (X \xrightarrow{id} X) \rightarrow \mathbb{B}^{i} (X \xrightarrow{id} X)$. \\
  
  We can view a bivariant theory as a way to codify and unify two theories, one covariant and one contravariant, in a cohesive and useful package. We will use this formalism when looking for covariant or contravariant versions of algebraic bordism (and cobordism).

\section{Algebraic Bordism and Cobordism}

We start by presenting Levine-Morel's definition of algebraic bordism applied to quasi-projective varieties over a base field $k$. This is the special case that will interest us in order to get a universal bivariant theory for quasi-projective varieties in the next section, and a bivariant version of algebraic bordism later. 

Let $X$ be a quasi-projective variety over the field $k$. A \textit{cobordism cycle} for $X$ is given by $(V \xrightarrow{h} X; L_1, \cdots, L_r)$, where $V$ is a smooth variety, integral over $k$, $h$ is a projective morphism, and $L_i$ are line bundles over $V$ (note that $r$ might be $0$). The dimension of one such cobordism cycle is defined as $\mathrm{dim}_k V - r$. Two cobordism cycles for $X$, $(V \xrightarrow{h} X; L_1, \cdots, L_r)$ and  $(V' \xrightarrow{h'} X; L'_1, \cdots, L'_r)$, are said to be \textit{equivalent} if there exists an isomorphism $g : V \rightarrow V'$ satisfying $h = h' \circ g$ and such that the pullbacks $g^*L'_1, \cdots, g^*L'_r$  give line bundles equivalent to $L_1, \cdots, L_r$ (but not necessarily in the same order). Denote an equivalence class of cobordism cycles by $[V \xrightarrow{h} X; L_1, \cdots, L_r]$.


Let $\mathcal{Z}_*(X)$ denote the free abelian group on the set of equivalence classes of cobordism cycles for $X$. This is a graded group, the \textit{group of cobordism cycles} of $X$.

Forgetting all about line bundles, we get the so-called \textit{cycle group} of $X$, $\mathcal{M}_*(X)$. This is the free abelian group on the set 
$$ \{[h : V \rightarrow X] : V {\rm \ is \ smooth \ and \ } h {\rm \ is \ projective} \},$$
graded by the dimension of $V$.

For a smooth variety $V$, its \textit{fundamental class} is defined by $[V] := [V \xrightarrow{c} pt]$, where $c$ is a smooth constant morphism to a point $pt$.

The algebraic bordism group for a quasi-projective $X$ will be constructed as a quotient of $\mathcal{Z}_* (X)$ by three relations of geometric flavour.

First, each line bundle $L$ over $X$ defines a \textit{first Chern class operator} (and thus an orientation) on $\mathcal{Z}_* (X)$ as a degree-one map
$\widetilde{c_1} (L) : \mathcal{Z}_* (X) \rightarrow \mathcal{Z}_{*-1} (X)$ by
$$\widetilde{c_1} (L) ([V \xrightarrow{h} X; L_1, \cdots, L_r]) = [V \xrightarrow{h} X; L_1, \cdots, L_r, h^*L]$$ 

Let $F_{\mathbb{L}}$ denote the universal formal group law of the Lazard ring $\mathbb{L}$ (see \cite{l} and \cite{q} for details).

\begin{axiom}(Dimension)
If $V$ is a smooth variety, integral over $k$, and $L_1, \cdots, L_r$ are line bundles over $V$ with $r > \mathrm{dim}_k V$, then
$\widetilde{c_1} (L_1) \circ \widetilde{c_1} (L_2) \circ \cdots \circ \widetilde{c_1} (L_r) ([V]) = 0$ in $\mathcal{Z}_* (V)$.
\end{axiom}

\begin{axiom}(Section)
If $V$ is a smooth variety, integral over $k$, $L$ is a line bundle over $V$ and $s : V \rightarrow L$ is a section transverse to the zero section of $L$, then
$\widetilde{c_1} (L) ([V]) = (i_Z)_* ([Z])$ in $\mathcal{Z}_* (V)$,
where $Z$ is defined as $s^{-1}(0)$ and $i_Z : Z \rightarrow V$ is the inclusion (which is a closed immersion).
\end{axiom}

\begin{axiom}(Formal Group Law)
If $V$ is a smooth variety, integral over $k$ and $L_1$ and $L_2$ are line bundles over $V$, then
$\widetilde{c_1} (L_1 \otimes L_2) = F_{\mathbb{L}}(\widetilde{c_1}(L_1),\widetilde{c_1}(L_2))$.
\end{axiom}

Suppose, for each $X$, that $\mathcal{R}_*(X) \subset \mathcal{Z}_* (X)$ is a subset formed by homogeneous elements. Then we can define the quotient $\mathcal{Z}_*/\mathcal{R}_*$, a functor on quasi-projective varieties, as follows. Let $\left<\mathcal{R}_*\right>(X) \subset \mathcal{Z}_* (X)$ be the subgroup generated by the set
$$\{ f_* \circ \widetilde{c_1} (L_1) \circ \cdots \widetilde{c_1} (L_r) \circ g^*(\rho) \}$$
where $Y$ and $Z$ are quasi-projective, $f:Y \rightarrow X$ is projective, $L_i$ are line bundles over $Y$, $g:Y \rightarrow Z$ is smooth and equidimensional, and $\rho \in \mathcal{R}_*(Z)$.

Define $\mathcal{Z}_*/\mathcal{R}_* (X)$ as the quotient group $\mathcal{Z}_*(X)/\left<\mathcal{R}_*\right>(X)$.

We say that the functor $\mathcal{Z}_*/\mathcal{R}_*$ is obtained from $\mathcal{Z}_*$ by \textit{imposing the relations} $\mathcal{R}_*(X)$ for all $X$. This construction can be made for any \textit{oriented Borel-Moore functor} (\cite{lm}, 2.1.3). This definition allows for the construction of algebraic bordism. This is done sequentially, by modding out each of the previous axioms.\\

$\mathcal{R}_*^\mathrm{Dim}(X) \subset \mathcal{Z}_*(X)$, for $X$ smooth, is defined as the subset of elements of the form $[V \xrightarrow{h} X; L_1, \cdots, L_r]$, with $\mathrm{dim}(V) < r$. Put $\underline{\mathcal{Z}}_*(X) : = \mathcal{Z}_*/\mathcal{R}_*^\mathrm{Dim}(X)$. \\

$\mathcal{R}_*^\mathrm{Sect}(X) \subset \underline{\mathcal{Z}}_*(X)$ is defined as the subset of elements of the form $$[V \xrightarrow{h} X; L_1, \cdots, L_r, L] - [Z \xrightarrow{hi} X; i^*L_1, \cdots, i^*L_r],$$ where $L$ is a line bundle over $V$, $s : V \rightarrow L$ is a section transverse to the zero section and $i : Z \rightarrow V$ is the closed immersion of the subvariety of zeros of $s$. Define the \textit{algebraic pre-cobordism group} $\underline{\Omega}_*(X)$ by $\underline{\mathcal{Z}}_*/\mathcal{R}_*^\mathrm{Sect}(X)$. \\

$\mathcal{R}_*^\mathrm{FGL}(X) \subset \mathbb{L}_* \, \otimes_\mathbb{Z} \, \underline{\Omega}_*(X)$ is defined as the subset of elements of the form 
$$(F_\mathbb{L_*}(\widetilde{c_1} (L_1),\widetilde{c_1} (L_2)) - \widetilde{c_1} (L_1 \otimes L_2))[\alpha],$$ where $X$ is smooth, $L_1$ and $L_2$ are line bundles over $X$ and $[\alpha] \in \underline{\Omega}_*(X)$. $\mathbb{L}\mathcal{R}_*^\mathrm{FGL}(X) \subset \mathbb{L}_* \otimes \underline{\Omega}_*(X)$ is the subset of elements of the form $a\beta$, with $a \in \mathbb{L}_*$ and $\beta \in \mathcal{R}_*^\mathrm{FGL}(X)$. Finally, the algebraic bordism group is defined as $\Omega_*(X) = \mathbb{L}_* \otimes \underline{\Omega}_*/\mathbb{L}\mathcal{R}_*^\mathrm{FGL}(X)$. \\

$\Omega_*$ can be viewed within the context of a more general concept, that of \textit{Borel-Moore functors}, where it can play a fundamental part as universal object. $\Omega_*$ gives an \textit{oriented Borel-Moore homology theory} (\cite{lm}, 5.1.3), which is universal for quasi-projective varieties over a field admitting resolution of singularities (\cite{lm}, 7.1.3).

If $f : X \rightarrow Y$ is a regular embedding (of codimension $d$) of smooth varieties over a field that admits resolution of singularities, we can get Gysin morphisms $f^* : \Omega_*(Y) \rightarrow \Omega_{* - d} (X)$ (\cite{lm}, 6.5.3). These allow to define a contravariant functor on smooth varieties (with regular embeddings as morphisms), algebraic cobordism, given by
$$\Omega^n (X) := \Omega_{n - \mathrm{dim } X} (X)$$
where, for $f : X \rightarrow Y$, the induced map $\Omega^* (Y) \rightarrow \Omega^* (X)$ is the corresponding Gysin map.

$\Omega^* (X)$ has a structure of graded ring (this comes from additional properties of $\Omega_*$), and defines what is called an \textit{oriented cohomology theory}. In fact, it is the universal such theory on the category of smooth varieties (provided one has resolution of singularities) (\cite{lm}, 7.1.3).

\section{A Universal Bivariant Theory for Quasi-projective Varieties}

Yokura and Sch\"urman define a universal bivariant theory in the category $\mathcal{V}$ of quasi-projective varieties over a base field $k$ (with no restriction on characteristic) (\cite{y} \cite{sy}).
We review their construction here, in order to apply it to a bivariant version of algebraic bordism. 

The class of confined maps will be that of proper morphisms, and the confined squares will be precisely all fiber squares.

For each morphism $X \xrightarrow{f} Y$ in $\mathcal{V}$, define ${\mathbb M}(\mathcal{V}/X \xrightarrow{f} Y)$ as the free abelian group generated by the set 
$$\{[h : W \rightarrow X] : W \in \mathcal{V}, h {\rm \ is \ proper \ and \ } f \circ h {\rm \ is \ smooth} \}$$

Here, $\left[ - \right]$ means equivalence class:  $h : W \rightarrow X$ and $h' : W' \rightarrow X$ are \textit{equivalent} if there exists an isomorphism $g : W \rightarrow W'$ such that $f \circ h' \circ g = f \circ h$. This group may be graded using the dimension of $W$.

Define the three operations of the bivariant theory as follows:

\begin{itemize}
  \item For the product \ $\bullet : {\mathbb M}(\mathcal{V}/X \xrightarrow{f} Y) \otimes {\mathbb M}(\mathcal{V}/Y \xrightarrow{g} Z) \rightarrow {\mathbb M}(\mathcal{V}/X \xrightarrow{gf} Z)$, consider $[h : W \rightarrow X] \in  {\mathbb M}(\mathcal{V}/X \xrightarrow{f} Y)$ and $[j : V \rightarrow Y] \in {\mathbb M}(\mathcal{V}/Y \xrightarrow{g} Z)$, and the diagram
  
  $$\xymatrix{
  & & V \ar[d]^{j} & \\
  W \ar[r]^{h} & X \ar[r]^{f} & Y \ar[r]^{g} & Z} $$ \\
Complete the diagram by adding the corresponding fibre squares, in order to obtain

 $$\xymatrix{
  W' \ar[r]^{h'} \ar[d]_{j''} & X' \ar[r]^{f'} \ar[d]_{j'} & V \ar[d]^{j} & \\
   W \ar[r]^{h} & X \ar[r]^{f} & Y \ar[r]^{g} & Z} $$ \\
Define then $[h : W \rightarrow X] \bullet [j : V \rightarrow Y]$ by $[hj'' : W' \rightarrow X]$, and extend bilinearly. \\

  \item For the pushforward $f_* : {\mathbb M}(\mathcal{V}/X \xrightarrow{gf} Z) \rightarrow {\mathbb M}(\mathcal{V}/Y \xrightarrow{g} Z)$, put $f_*([h : W \rightarrow X]) = [fh : W \rightarrow Y]$, and extend linearly. \\ 
  \item For the pullback $g^* : {\mathbb M}(\mathcal{V}/X' \xrightarrow{h'} Y') \rightarrow {\mathbb M}(\mathcal{V}/X \xrightarrow{h} Y)$, given a confined square
  $\xymatrix{
  X \ar[r]^{f} \ar[d]_{h} & X' \ar[d]^{h'} \\
  Y \ar[r]^{g}  & Y'  }$, \\
  consider $[j : W' \rightarrow X'] \in {\mathbb M}(\mathcal{V}/X' \xrightarrow{h'} Y')$, and the diagram
  
 $\xymatrix{
 & W' \ar[d]^{j} \\
  X \ar[r]^{f} \ar[d]_{h} & X' \ar[d]^{h'} \\
  Y \ar[r]^{g}  & Y'  }$, \\ \\
Complete the diagram by filling in the top fibre square, in order to obtain

 $\xymatrix{
 W \ar[r]^{f'} \ar[d]_{j'} & W' \ar[d]^{j} \\
  X \ar[r]^{f} \ar[d]_{h} & X' \ar[d]^{h'} \\
  Y \ar[r]^{g}  & Y'  }$, \\
  
Put then $g^*([j : W' \rightarrow X']) = [j' : W \rightarrow X]$, and extend linearly. \\ 
\end{itemize}
 
\begin{theorem} (\cite{sy} \cite{y})
${\mathbb M}(\mathcal{V}/X \xrightarrow{f} Y)$ determines a bivariant theory.
\end{theorem}


Yokura and Sch\"urman call ${\mathbb M}(\mathcal{V}/X \xrightarrow{f} Y)$ the {\it pre-motivic bivariant Grothendieck group} for $f : X \rightarrow Y$. This is a bivariant version of the previous cycle group $\mathcal{M}_*(X)$. 

This bivariant theory is universal among all similar bifunctors defined on the category of quasi-projective varieties over a fixed field $k$ (\cite{sy}, Thm 3.3). \\ 

A \textit{bivariant cobordism cycle} for $f : X \rightarrow Y$ is given by $(V \xrightarrow{h} X; L_1, \cdots, L_r)$, where $V$ is a smooth variety, $h$ is a projective morphism, and $L_i$ are line bundles over $V$. The dimension of one such cobordism cycle is defined as $\mathrm{dim} V - r$. Two cobordism cycles for $f : X \rightarrow Y$, $(V \xrightarrow{h} X; L_1, \cdots, L_r)$ and  $(V' \xrightarrow{h'} X; L'_1, \cdots, L'_r)$, are said to be \textit{equivalent} if there exists an isomorphism $g : V \rightarrow V'$ satisfying $f \circ h = f \circ h' \circ g$ and such that the pullbacks $g^*L'_1, \cdots, g^*L'_k$  give line bundles equivalent to $L_1, \cdots, L_k$ (but not necessarily in the same order). Denote an equivalence class of cobordism cycles by $[V \xrightarrow{h} X; L_1, \cdots, L_k]_f$. The dependence on $f$ is explicit in the notation.

Define $\mathbb{OB}(X \xrightarrow{f} Y)$ as the free abelian group generated by the set 
$$\{[h : W \rightarrow X; L_1, \cdots, L_r]_f : W \in \mathcal{V}, h {\rm \ projective, \ } f \circ h {\rm \ smooth, \ } L_i {\ \rm line \ bundles \ over \ } W  \}$$

$\mathbb{OB}(X \xrightarrow{f} Y)$ is a bivariant analogue of the previous group of cobordism cycles $\mathcal{Z}_*(X)$ from which algebraic bordism was defined and, as such, one can naturally define a (bivariant) orientation in this context. For each line bundle $L$ over $X$, the first Chern class operator
$$\widetilde{c_1} (L) : \mathbb{OB}(X \xrightarrow{f} Y) \rightarrow \mathbb{OB} (X \xrightarrow{f} Y)$$
will be given by 
$$\widetilde{c_1} (L) ([V \xrightarrow{h} X; L_1, \cdots, L_r]_f) = [V \xrightarrow{h} X; L_1, \cdots, L_r, h^*L]_f$$ 

With this orientation, $\mathbb{OB}$ becomes a universal oriented bivariant theory (\cite{y}, Thm. 3.1). If we define the degree of $[h : W \rightarrow X; L_1, \cdots, L_r] \in \mathbb{OB} (X \xrightarrow{f} Y)$ by $\rm{dim \ } W - r$, then each $\widetilde{c_1} (L) $ is a degree-one map.

The covariant theory $\mathbb{OB}_*$ associated to $\mathbb{OB}$ is $\mathbb{OB}_*(X) = \mathbb{OB}(X \xrightarrow{c} pt) \cong \mathcal{Z}_*(X)$, the previous group of cobordism cycles.

As for the contravariant theory $\mathbb{OB}^*$, each $\mathbb{OB}^*(X) = \mathbb{OB}(X \xrightarrow{id} X)$ is generated from elements $[V \xrightarrow{h} X]_{id}$, with $h$ smooth and projective, by application of Chern class operators $\widetilde{c_1} (L)$ for all line bundles $L$ over $X$.

\section{A Contravariant Version of Algebraic Bordism}

To generalize the three axioms present in the definition of algebraic bordism, we start by extending the definition of  fundamental class for smooth varieties.

Let $X$ be smooth and $X \xrightarrow{f} Y$ be a morphism in $\mathcal{V}$. Define the fundamental class of $X \xrightarrow{f} Y$ by $[X \xrightarrow{f} Y] := [X \xrightarrow{id} X]_f \in \mathbb{OB}(X \xrightarrow{f} Y)$, the equivalence class of the bivariant cobordism cycle $(X \xrightarrow{id} X)$ (no line bundles included). Its dimension is of course the dimension of $X$. This definition agrees with the previous concept: if $X$ is smooth, $[X]$ coincides with $[X \xrightarrow{c} pt]$ by construction.
  
Unraveling the definition further, we get 
$$[X \xrightarrow{f} Y] = [X \xrightarrow{id} X]_f = \{ (V \xrightarrow{h} X) : \exists \textrm{\ isomorphism \ } g : V \rightarrow X \textrm{\ with } f \circ h = f \circ g \}.$$

The dependence of the fundamental class $[X \xrightarrow{f} Y]$ on $Y$ that comes from the definition of equivalent bivariant cobordism cycles is what permits a generalization of the axioms included in the definition of algebraic bordism to a universal bivariant context. These new axioms are presented below. 
 
\begin{axiom}(Dimension)
If $X$ is a smooth variety, $X \xrightarrow{f} Y$ is in $\mathcal{V}$, and $L_1, \cdots, L_r$ are line bundles over $X$ with $r > \mathrm{dim} X$, then
$\widetilde{c_1} (L_1) \circ \widetilde{c_1} (L_2) \circ \cdots \circ \widetilde{c_1} (L_r) ([X \xrightarrow{f} Y]) = 0$ in $\mathbb{OB}(X \xrightarrow{f} Y)$.
\end{axiom} 

Let $X$ be smooth and $i_Z : Z \rightarrow X$ be a closed immersion of a subvariety. Since $i_Z$ is proper, we can define the pushforward $(i_Z)_* : \mathbb{OB}(Z \xrightarrow{fi_Z} Y) \rightarrow \mathbb{OB}(X \xrightarrow{f} Y)$.

\begin{axiom}(Section)
If $X$ is a smooth variety, $X \xrightarrow{f} Y$ is in $\mathcal{V}$, $L$ is a line bundle over $X$ and $s : X \rightarrow L$ is a section transverse to the zero section of $L$, \\ then
$\widetilde{c_1} (L) ([X \xrightarrow{f} Y]) = (i_Z)_* ([Z \xrightarrow{f i_Z} Y])$ in $\mathbb{OB}(X \xrightarrow{f} Y)$, \\
where $Z$ is defined as $s^{-1}(0)$ and $i_Z : Z \rightarrow X$ is the inclusion (which is a closed immersion).
\end{axiom}

\begin{axiom}(Formal Group Law)
If $X$ is a smooth variety,  $X \xrightarrow{f} Y$ is in $\mathcal{V}$, and $L_1$ and $L_2$ are line bundles over $X$, then
$\widetilde{c_1} (L_1 \otimes L_2) = F_{\mathbb{L}}(\widetilde{c_1}(L_1),\widetilde{c_1}(L_2))$ as operators $\mathbb{L}_* \otimes_{\mathbb{Z}} \mathbb{OB}(X \xrightarrow{f} Y) \rightarrow \mathbb{L}_* \otimes_{\mathbb{Z}} \mathbb{OB}(X \xrightarrow{f} Y)$,
where $F_{\mathbb{L}}$ denotes the universal formal group law of the Lazard ring $\mathbb{L}_*$.
\end{axiom}

The Chern class operator $\widetilde{c_1} (L_1 \otimes L_2)$ in the last axiom should be considered as acting on $\mathbb{L}_* \otimes_{\mathbb{Z}} \mathbb{OB}(X \xrightarrow{f} Y)$ (the action is trivial on $\mathbb{L}_*$). Note that $F_{\mathbb{L}}(\widetilde{c_1}(L_1),\widetilde{c_1}(L_2))$ has coefficients in $\mathbb{L}_*$.

From these axioms, we are now ready to provide a bivariant analogue of algebraic bordism. We will be applying each axiom in succession, and so the result will be a quotient of $\mathbb{OB}$.

Suppose, for each $f : X \rightarrow Y$ in $\mathcal{V}$, that $\mathcal{R}_*(X \xrightarrow{f} Y) \subset \mathbb{OB}(X \xrightarrow{f} Y)$ is a subset formed by homogeneous elements. We can define the quotient $\mathbb{OB}/\mathcal{R}_*$, an oriented bivariant theory on $\mathcal{V}$, as follows. Let $\left<\mathcal{R}_*\right>(X \xrightarrow{f} Y) \subset\mathbb{OB}(X \xrightarrow{f} Y)$ be the subgroup generated by the set
$$\{ f'_* \circ \widetilde{c_1} (L_1) \circ \cdots \widetilde{c_1} (L_r) \circ g^*(\rho) \}$$
which relates to

$\xymatrix{
 X' \ar[r]^{} \ar[d]_{f'} & Z' \ar[dd]^{z} \\
  X \ar[d]_{f} & \\
  Y  \ar[r]^{g} & Z }$ \\

where $f'$ and the outside square are confined, $f:Y \rightarrow X$ is projective, $L_i$ are line bundles over $Y$, $g:Y \rightarrow Z$ is smooth and equidimensional, and $\rho \in \mathcal{R}_*(Z' \xrightarrow{z} Z)$.

Define $\mathbb{OB}/\mathcal{R}_* (X \xrightarrow{f} Y)$ as the quotient group $\mathbb{OB}(X \xrightarrow{f} Y)/\left<\mathcal{R}_*\right>(X \xrightarrow{f} Y)$.

As before, we say that $\mathbb{OB}/\mathcal{R}_*$ is obtained from $\mathbb{OB}$ by \textit{imposing the relations} $\mathcal{R}_*(X \xrightarrow{f} Y)$ for all $X \xrightarrow{f} Y$. This construction can be made for any oriented bivariant theory.\\

$\mathcal{R}_*^\mathrm{Dim}(X \xrightarrow{f} Y) \subset \mathbb{OB}(X \xrightarrow{f} Y)$, for $X$ smooth, is defined as the subset of elements of the form $\widetilde{c_1} (L_1) \circ \widetilde{c_1} (L_2) \circ \cdots \circ \widetilde{c_1} (L_r) ([X \xrightarrow{f} Y])$, with $\mathrm{dim}(X) < r$. Put $\underline{\mathcal{Z}}_*(X \xrightarrow{f} Y) : = \mathbb{OB}/\mathcal{R}_*^\mathrm{Dim}(X \xrightarrow{f} Y)$. \\

$\mathcal{R}_*^\mathrm{Sect}(X \xrightarrow{f} Y) \subset \underline{\mathcal{Z}}_*(X \xrightarrow{f} Y)$ is defined as the subset of elements of the form $\widetilde{c_1} (L) ([X \xrightarrow{f} Y]) - (i_Z)_* ([Z \xrightarrow{f i_Z} Y])$, where $X$ is smooth, $L$ is a line bundle over $X$, $s : X \rightarrow L$ is a section transverse to the zero section of $L$, and $i _Z : Z \rightarrow X$ is the closed immersion of the subvariety of zeros of $s$. View $(i_Z)_* ([Z \xrightarrow{f i_Z} Y])$ as $[Z \xrightarrow{id} Z]_{fi_Z} = [i_Z(Z) \xrightarrow{id} i_Z(Z)]_{f}$. Define $\underline{\Omega}_*(X \xrightarrow{f} Y)$ by $\underline{\mathcal{Z}}_*/\mathcal{R}_*^\mathrm{Sect}(X \xrightarrow{f} Y)$. \\

$\mathcal{R}_*^\mathrm{FGL}(X \xrightarrow{f} Y) \subset \mathbb{L}_* \, \otimes_\mathbb{Z} \, \underline{\Omega}_*(X \xrightarrow{f} Y)$ is defined as the subset of elements of the form 
$(F_\mathbb{L_*}(\widetilde{c_1} (L_1),\widetilde{c_1} (L_2)) - \widetilde{c_1} (L_1 \otimes L_2))[\alpha]$, where $X$ is smooth, $L_1$ and $L_2$ are line bundles over $X$ and $[\alpha] \in \underline{\Omega}_*(X \xrightarrow{f} Y)$. $\mathbb{L}\mathcal{R}_*^\mathrm{FGL}(X \xrightarrow{f} Y) \subset \mathbb{L}_* \otimes \underline{\Omega}_*(X \xrightarrow{f} Y)$ is the subset of elements of the form $a\beta$, with $a \in \mathbb{L}_*$ and $\beta \in \mathcal{R}_*^\mathrm{FGL}(X \xrightarrow{f} Y)$. \\

Our first bivariant theory generalizing algebraic bordism will then be defined as \\ $\mathbb{OB}1(X \xrightarrow{f} Y) := \mathbb{L}_* \otimes \underline{\Omega}_*/\mathbb{L}\mathcal{R}_*^\mathrm{FGL}(X \xrightarrow{f} Y)$. \\

It is immediate to notice that $\mathbb{OB}1_*(X) := \mathbb{OB}1(X \xrightarrow{f} pt)$ is a covariant theory that gives a group isomorphic to $\Omega_*(X)$ for each quasi-projective variety. \\

We now describe the associated contravariant theory $\mathbb{OB}1^*(X) := \mathbb{OB}1(X \xrightarrow{id} X)$ for smooth varieties. If $X$ is smooth, we saw that $\mathbb{OB}^*(X) := \mathbb{OB}(X \xrightarrow{id} X)$ is given by 
$$\{ [V \xrightarrow{h} X; L_1, \cdots L_r ] : W \in \mathcal{V}, h {\rm \ smooth \ and \ projective, \ } L_i {\ \rm line \ bundles \ over \ } W \}^+.$$ 

In order to get $\mathbb{OB}1^*(X)$, we will mod this group by $\mathcal{R}_*(X \xrightarrow{id} X)$ for each of the three sets $\mathcal{R}_*$ (one for each axiom). \\

The first axiom requires that $r$ must always be at most $\rm{dim \ }X$ in each cobordism cycle. \\

For the second axiom, we calculate first
$$[X \xrightarrow{id} X] = [X \xrightarrow{id} X]_{id} = \{ (V \xrightarrow{h} X) : \exists \textrm{\ isomorphism \ } g : V \rightarrow X \textrm{\ with \ } h = g \}$$

That is,  $[X \xrightarrow{id} X] = \{ (V \xrightarrow{h} X) : V \textrm{\ smooth}, h \textrm{\ isomorphism}\}.$ \\

Also, for $Z$ as was considered in $\mathcal{R}_*^\mathrm{Sect}$, 
$$[Z \xrightarrow{i_Z} X] = [Z \xrightarrow{id} Z]_{i_Z} = \{ (V \xrightarrow{h} Z) : \exists \textrm{\ isomorphism \ } g : V \rightarrow Z \textrm{\ with \ } i_Z \circ h = i_Z \circ g \}$$  

and $(i_Z)_* ([Z \xrightarrow{i_Z} X]) = [i_Z(Z) \xrightarrow{id} i_Z(Z)]_{id}$. \\

The bivariant pre-cobordism group in this case will be defined from $\mathbb{OB}$ by considering bivariant cobordism cycles where $r$ must always be at most $\rm{dim \ }X$ and such that, for all $Z$ and $L$ under the mentioned conditions, we get $[X \xrightarrow{id} X; L]_{id} = [i_Z(Z) \xrightarrow{id} i_Z(Z)]_{id}$. \\

For $\mathbb{OB}1_*(X)$, we additionally impose (formally) the universal formal group law. 

\section{A Covariant Version of Algebraic Cobordism}

In certain cases, we can define the dual bivariant theory to a given (graded) bivariant theory. For our purposes, this allows the definition of a bivariant theory for quasi-projective varieties ${\mathbb{M}}'(\mathcal{V}/X \xrightarrow{f} Y)$ from which we get another bivariant theory $\mathbb{OB}2(X \xrightarrow{f} Y)$ whenever $X$ is smooth, which satisfies $\mathbb{OB}2(X \xrightarrow{id} X) \cong \Omega^*(X)$.

Start by defining ${\mathbb M}'(\mathcal{V}/X \xrightarrow{f} Y)$ as the free abelian group generated by the set 
$$\{[h : Y \rightarrow W] : W \in \mathcal{V}, h {\rm \ is \ proper \ and \ } h \circ f {\rm \ is \ smooth} \},$$

(with same definition of equivalence class as before).

We can put on ${\mathbb M}'$ the structure of a bivariant theory. For confined maps, take those $f : X \rightarrow Y$ such that, for any $g : Y \rightarrow Z$, $gf$ smooth implies $g$ smooth. The confined squares will be all fiber squares. Define the three required operations as follows:

\begin{itemize}
  \item For the product \ $\bullet : {\mathbb M}'(\mathcal{V}/X \xrightarrow{f} Y) \otimes {\mathbb M}'(\mathcal{V}/Y \xrightarrow{g} Z) \rightarrow {\mathbb M}'(\mathcal{V}/X \xrightarrow{gf} Z)$, consider $[h : Y \rightarrow W] \in  {\mathbb M}'(\mathcal{V}/X \xrightarrow{f} Y)$ and $[j : Z \rightarrow V] \in {\mathbb M}'(\mathcal{V}/Y \xrightarrow{g} Z)$, and the diagram
  
  $$\xymatrix{
   X \ar[r]^{f} & Y \ar[r]^{g} \ar[d]^{h} & Z \ar[r]^{j} & V \\
   & W & &} \\
  $$ \\
Complete the diagram by adding the corresponding pushout squares, in order to obtain

 $$\xymatrix{
  X \ar[r]^{f} & Y \ar[r]^{g} \ar[d]^{h} & Z \ar[r]^{j} \ar[d]^{h'} & \ar[d]^{h''} V \\
     & W \ar[r]^{g'} & Z' \ar[r]^{j'} & V'} $$ \\
Define then $[h : Y \rightarrow W] \bullet [j : Z \rightarrow V]$ by $[j'h' : Z \rightarrow V']$, and extend bilinearly. \\

  \item For the pushforward $f_* : {\mathbb M}'(\mathcal{V}/X \xrightarrow{gf} Z) \rightarrow {\mathbb M}'(\mathcal{V}/Y \xrightarrow{g} Z)$, put $f_*([h : Z \rightarrow W]) : = [h : Z \rightarrow W]$, and extend linearly. Note that $f$, being confined, implies that $hg$ is smooth. $f_*$ should not be viewed as an identity.\\ 
  
  \item For the pullback $g^* : {\mathbb M}'(\mathcal{V}/X' \xrightarrow{h'} Y') \rightarrow {\mathbb M}'(\mathcal{V}/X \xrightarrow{h} Y)$, given a confined square
  $\xymatrix{
  X \ar[r]^{f} \ar[d]_{h} & X' \ar[d]^{h'} \\
  Y \ar[r]^{g}  & Y'  }$, \\
  consider $[j : Y' \rightarrow W] \in {\mathbb M}'(\mathcal{V}/X' \xrightarrow{h'} Y')$, and the diagram
  
 $\xymatrix{
  X \ar[r]^{f} \ar[d]_{h} & X' \ar[d]^{h'} \\
  Y \ar[r]^{g}  & Y' \ar[d]_{j} \\
  & W }$, \\
  
Put then $g^*([j : Y' \rightarrow W]) : = [jg : Y \rightarrow W]$, and extend linearly. \\
\end{itemize}
 
It can be seen, similarly to what happens with ${\mathbb M}$, that these operations satisfy the required axioms, and that ${\mathbb M}'(\mathcal{V}/X \xrightarrow{f} Y)$ in fact determines a bivariant theory. \\

If $X$ is smooth, each $[h : Y \rightarrow W] \in {\mathbb M}'(\mathcal{V}/X \xrightarrow{f} Y)$ can actually be considered as an element in ${\mathbb M}(\mathcal{V}/Y \xrightarrow{h} W)$. The image of this element  in the quotient by the previous axioms belongs then to $\mathbb{OB}1 (Y \xrightarrow{h} W)$. Call it $I([h : Y \rightarrow W])$. 

We define, for each $X \xrightarrow{f} Y$ in $\mathcal{V}$ such that $X$ is smooth, a new bivariant theory
$$\mathbb{OB}2(X \xrightarrow{f} Y) := \{ I([h : Y \rightarrow W]) \in \mathbb{OB}1 (Y \xrightarrow{h} W) : [h : Y \rightarrow W] \in {\mathbb M}'(\mathcal{V}/X \xrightarrow{f} Y) \}^+$$

The axioms used for our previous bivariant algebraic bordism are considered in this definition but now, contrary to what happened with $\mathbb{OB}1$, $\mathbb{OB}2$ is not directly seen as a quotient of ${\mathbb M}'$. It can in a sense be considered as a quotient of $\mathbb{OB}$, as we'll see in section 10. \\

If $X$ is smooth, we get also
$$\mathbb{M}'(\mathcal{V}/X \xrightarrow{id} X) = \{[h : X \rightarrow W] : W \in \mathcal{V}, h {\rm \ is \ proper \ and \ smooth} \}$$

This defines a contravariant functor (by composing on the left), and we can check that $\mathbb{OB}2(X \xrightarrow{id} X)$ (graded by the dimension of $W$) is isomorphic to the algebraic cobordism (ring)  $\Omega^*(X)$. To see this, take a cycle $[h : X \rightarrow W] \in \mathbb{M}'(\mathcal{V}/X \xrightarrow{id} X)$ (of dimension $n = \mathrm{dim } W$). Note that, since $X$ is smooth, $h$ defines also a cycle in $\Omega_{\mathrm{dim } X}(W)$. Using the Gysin map $h^* : \Omega_* (W) \rightarrow \Omega_{*-d}(X)$, where $d = \mathrm{dim } W - \mathrm{dim } X$, we can obtain an element $h^*([h : X \rightarrow W]) \in \Omega_{\mathrm{dim } W - d}(X) = \Omega_{\mathrm{dim } X}(X)$. The correspondence $\{ I([h : X \rightarrow W])\} \mapsto h^*([h : X \rightarrow W])$ gives an isomorphism $\mathbb{OB}2(X \xrightarrow{id} X) \rightarrow \Omega^*(X)$. \\

Since algebraic cobordism was not the contravariant version of the bivariant theory $\mathbb{OB}1$ (when applied to smooth varieties) in the previous section, the covariant version of $\mathbb{OB}2$, given by $\mathbb{OB}2(X \xrightarrow{c} pt)$, for $X$ smooth, will not be algebraic bordism.

We describe this new covariant theory for smooth varieties: 
$${\mathbb M}'(\mathcal{V}/X \xrightarrow{c} pt) = \{[h : pt \rightarrow W] : W \in \mathcal{V}, h {\rm \ is \ proper \ and \ } h \circ c {\rm \ is \ smooth} \}^+.$$ 

That is, ${\mathbb M}'(\mathcal{V}/X \xrightarrow{c} pt)= \{ h \in W : W \in \mathcal{V}\}^+$, the free group on all points of all quasi-projective varieties of $\mathcal{V}$. \\

This means that we get also

$$\mathbb{OB}2(X \xrightarrow{c} pt) = \{ h \in W : W \in \mathcal{V}\}^+,$$ 
which does not depend on the original $X$. We can define thus this covariant theory also when $X$ is any quasi-projective variety.

\section{Double Point Bordism and Cobordism}

Levine and Pandharipande provide in \cite{lp} a new construction that results, for any variety $X$, in a group isomorphic to the algebraic bordism group $\Omega_* (X)$. This new group is called \textit{double point bordism} and has a deeper geometric description than the more formal one coming from $\Omega_*$. It is done from $\mathcal{M}_*(X)$ by modding out the \textit{double point relations} over $X$.

Recall the definition of the cycle group for any variety $X$: 
$$\mathcal{M}_*(X) = \{[h : V \rightarrow X] : V {\rm \ is \ smooth \ and \ } h {\rm \ is \ projective} \}^+ \\$$ 

It is graded by the dimension of $V$. \\

Now work over a field of characteristic zero. Let $Y$ be a smooth variety. A \textit{double point degeneration} over $0 \in \mathbb{P}^1$ is a morphism $\pi : Y \rightarrow \mathbb{P}^1$ such that $\pi^{-1}(0) = A \cup B$ for two smooth subvarieties of $Y$ of codimension $1$ that intersect transversely. $D = A \cap B$ is the \textit{double point locus} of $\pi$ over $0$. Taking the normal bundles of $D$ in $A$ and $B$, $N_{A/D}$ and $N_{B/D}$, we get two isomorphic projective bundles $\mathbb{P}(\mathcal{O}_D \otimes N_{A/D}) \rightarrow D$ and  $\mathbb{P}(\mathcal{O}_D \otimes N_{B/D}) \rightarrow D$, both denoted by $P : \mathbb{P}(\mathbb{\pi}) \rightarrow D$. 

Let $Y$ be a smooth variety and $X$ any variety, and consider a projective morphism $g : Y \rightarrow X \times \mathbb{P}^1$ such that composition with projection onto the second factor gives a double point degeneration $\pi$ over $0 \in \mathbb{P}^1$. If $\xi \in \mathbb{P}^1$ is any regular value of $\pi$, denote the fiber over $\xi$ by $Y_\xi$. $g$ is called a \textit{double point cobordism} with degenerate fiber over $0$ and smooth fiber over $\xi$, and has its \textit{associated double point relation} over $X$: 

$$\left[Y_\xi \xrightarrow{p_1g} X\right] - \left[A \xrightarrow{p_1g} X\right] - \left[B \xrightarrow{p_1g} X\right] + \left[\mathbb{P}(\pi) \xrightarrow{p_1gP} X \right],$$ 

where $p_1$ is projection onto the first factor.

Denote by $\mathcal{R}_*(X)$ the subgroup of $M_*(X)$ generated by all the double point relations over $X$, and by $\omega_*(X)$ the quotient $M_*(X)/\mathcal{R}_*(X)$. $\omega_*(X)$ is called the \textit{double point cobordism group} for the variety $X$.

\begin{theorem} (\cite{lp}, Thm. 0.1)
For any quasi-projective variety $X$, $\Omega_*(X) \cong \omega_*(X)$.
\end{theorem}

$\omega_*$ admits Chern class operators and an orientation, whose definition is not as formal and direct as the one for $\Omega_*$ (\cite{lp}, Section 6).

For a \textit{smooth} variety $X$, we can, as before, define a contravariant functor $\omega^*$ that gives $\omega^*(X) = \omega_{* - \mathrm{dim } X}(X)$, which is called the \textit{double point cobordism ring} for $X$ (and which will of course be isomorphic to the algebraic cobordism ring $\Omega^*(X)$).

\section{A Contravariant Version of Double Point Bordism}

In this section, we generalize the double point relations to a bivariant setting, in order to get a contravariant analogue to double point bordism.

Recall that, for each morphism $X \xrightarrow{f} Y$ in $\mathcal{V}$, we defined the pre-motivic bivariant Grothendieck group
$${\mathbb M}(\mathcal{V}/X \xrightarrow{f} Y) = \{[h : W \rightarrow X] : W \in \mathcal{V}, h {\rm \ is \ proper \ and \ } f \circ h {\rm \ is \ smooth} \}^+.$$

We will get a bivariant theory $\mathbb{OB}3$ once we form a quotient of ${\mathbb M}(\mathcal{V}/X \xrightarrow{f} Y)$ by the corresponding double point relations.

Let $V$ be any quasi-projective variety. Consider a morphism $\pi : V \rightarrow \mathbb{P}^1$ such that $\pi^{-1}(0) = A \cup B$ for two subvarieties of $V$ of codimension $1$ that intersect transversely. As before, get the projective bundle $P : \mathbb{P}(\mathbb{\pi}) \rightarrow D$.

Fix a $X \xrightarrow{f} Y$ in $\mathcal{V}$. Let $[h : V \rightarrow X] \in {\mathbb M}(\mathcal{V}/X \xrightarrow{f} Y)$, and consider a projective morphism $g : V \rightarrow X \times \mathbb{P}^1$ such that the composition 
$V \xrightarrow{g} X \times \mathbb{P}^1 \xrightarrow{p_2}  \mathbb{P}^1$ gives a morphism $\pi$ with the properties presented in the previous paragraph, and such that the compositions $V_\xi \xrightarrow{fp_1g} Y$,  $A \xrightarrow{fp_1g} Y$, $B \xrightarrow{fp_1g} Y$, and $\mathbb{P}(\pi) \xrightarrow{fp_1gP} Y$ are all smooth. If $\xi \in \mathbb{P}^1$ is any regular value of $\pi$, denote the fiber over $\xi$ by $V_\xi$. We can call $g$ a \textit{double point cobordism relative to $f$}. Its \textit{associated double point relation} over $X \xrightarrow{f} Y$ is: 

$$\left[V_\xi \xrightarrow{p_1g} X\right] - \left[A \xrightarrow{p_1g} X\right] - \left[B \xrightarrow{p_1g} X\right] + \left[\mathbb{P}(\pi) \xrightarrow{p_1gP} X\right]$$

The elements present in such relations belong to ${\mathbb M}(\mathcal{V}/X \xrightarrow{f} Y)$.

As before, denote by $\mathcal{R}_*(X \xrightarrow{f} Y)$ the subgroup of ${\mathbb M}(\mathcal{V}/X \xrightarrow{f} Y)$ generated by all double point relations over $X \xrightarrow{f} Y$, and by $\mathbb{OB}3(X \xrightarrow{f} Y)$ the quotient ${\mathbb M}(X \xrightarrow{f} Y)/\mathcal{R}_*(X \xrightarrow{f} Y)$. 

If $Y$ is a point, a double point relation relative to $X \xrightarrow{c} pt$ is just a double point relation for $X$, and so $\mathcal{R}_*(X \xrightarrow{c} pt) \cong \mathcal{R}_*(X)$. This gives $$\mathbb{OB}3_*(X) := \mathbb{OB}3(X \xrightarrow{c} pt) = {\mathbb M}(X \xrightarrow{c} pt)/\mathcal{R}_*(X \xrightarrow{c} pt) \cong  \mathcal{M}_* (X)/\mathcal{R}_*(X) \cong \omega_*(X),$$ the double point bordism group of $X$.

For the contravariant theory, we analyze the double point relations over $X \xrightarrow{id} X$ whenever $X$ is a smooth variety. First, we have \\
$$ {\mathbb M}(X \xrightarrow{id} X) = \{ [h : W \rightarrow X] : W \in \mathcal{V}, h {\ \rm{proper \ and \ smooth}}\}^+.$$ \\
Notice that each element $[h : W \rightarrow X]$ in this group also belongs to $ {\mathbb M}'(W \xrightarrow{id} W)$. \\
 
Applying the definitions, we see that $\mathcal{R}_*(X \xrightarrow{id} X)$ is given by double point relations of the form \\
$$\left[V_\xi \xrightarrow{p_1g} X\right] - \left[A \xrightarrow{p_1g} X\right] - \left[B \xrightarrow{p_1g} X\right] + \left[\mathbb{P}(\pi) \xrightarrow{p_1gP} X\right],$$
with all four terms belonging to ${\mathbb M}(X \xrightarrow{id} X)$, that is, such that all maps featured in those terms are smooth (and proper). These are relations that are similar to those appearing in the definition of (total) double point bordism, but in this case the domain of $g : V \rightarrow X \times \mathbb{P}^1$ need not be smooth.

\section{A Covariant Version of Double Point Cobordism}

We can also obtain a covariant version of double point cobordism for smooth varieties, following a process similar to the one carried in section 6. This implies defining a new bivariant theory $\mathbb{OB}4(X \xrightarrow{f} Y)$ such that $\mathbb{OB}4(X \xrightarrow{c} X) \cong \omega^*(X)$ for $X$ smooth. The new covariant theory, not equivalent to double point bordism, is the covariant version $\mathbb{OB}4_*(X) := \mathbb{OB}4(X \xrightarrow{c} pt)$ of this new bivariant theory.

If $X$ is smooth, each $[h : Y \rightarrow W] \in {\mathbb M}'(\mathcal{V}/X \xrightarrow{f} Y)$ can actually be considered as an element in ${\mathbb M}(\mathcal{V}/Y \xrightarrow{h} W)$. The image of this element  in the quotient by the previous bivariant double point relations belongs then to $\mathbb{OB}3 (Y \xrightarrow{h} W)$. Call it $J([h : Y \rightarrow W])$. 

We define, for each $X \xrightarrow{f} Y$ in $\mathcal{V}$ such that $X$ is smooth, a new bivariant theory
$$\mathbb{OB}4(X \xrightarrow{f} Y) := \{ J([h : Y \rightarrow W]) \in \mathbb{OB}3 (Y \xrightarrow{h} W) : [h : Y \rightarrow W] \in {\mathbb M}'(\mathcal{V}/X \xrightarrow{f} Y) \}^+$$ 

The double point relations used for our previous bivariant double point bordism are considered in this definition but now, as for $\mathbb{OB}2$, $\mathbb{OB}4$ is not directly seen as a quotient of ${\mathbb M}'$. Again, it can in a sense be considered as a quotient of $\mathbb{OB}$, as we'll see in section 10. \\

If $X$ is smooth, we get also
$$\mathbb{M}'(\mathcal{V}/X \xrightarrow{id} X) = \{[h : X \rightarrow W] : W \in \mathcal{V}, h {\rm \ is \ proper \ and \ smooth} \}$$

This defines a contravariant functor, again by composing on the left, such that the associated $\mathbb{OB}4(X \xrightarrow{id} X)$ (graded by the dimension of $W$) is isomorphic to the algebraic double point cobordism ring  $\omega^*(X)$. Take a cycle $[h : X \rightarrow W] \in \mathbb{M}'(\mathcal{V}/X \xrightarrow{id} X)$ (of dimension $n = \mathrm{dim } W$). Using the similar reasonings of section $6$, we can obtain an element $h^*([h : X \rightarrow W]) \in \Omega_{\mathrm{dim } X}(X)$, and the isomorphism $\phi : \Omega_*(X) \cong \omega_*(X)$ finally gives an element $\phi h^*([h : X \rightarrow W]) \in \omega_{\mathrm{dim } X}(X)$. The correspondence $\{ J([h : X \rightarrow W])\} \mapsto \phi g^*([h : X \rightarrow W])$ is an isomorphism $\mathbb{OB}4(X \xrightarrow{id} X) \rightarrow \omega^*(X)$ (deducing the indexing from that of $\omega_*(X)$). \\

Since algebraic cobordism was not the contravariant version of the bivariant theory $\mathbb{OB}3$ (when applied to smooth varieties) in the previous section, the covariant version of $\mathbb{OB}4$, given by $\mathbb{OB}4(X \xrightarrow{c} pt)$, for $X$ smooth, will not be algebraic bordism.

For this new covariant theory $\mathbb{OB}4(X \xrightarrow{c} pt)$(with $X$ smooth), recall that 
$${\mathbb M}'(\mathcal{V}/X \xrightarrow{c} pt) = \{ h \in W : W \in \mathcal{V}\}^+,$$ 
the free group on all points of all quasi-projective varieties of $\mathcal{V}$. This directly gives
$$\mathbb{OB}4(X \xrightarrow{c} pt) = \{ h \in W : W \in \mathcal{V}\}^+,$$ 
which does not depend on the original $X$ and allows thus an extension of this covariant theory for any quasi-projective $X$.

\section{Relations Between the Presented Bivariant Theories}

We end with a brief account of how the presented bivariant theories and corresponding covariant and contravariant parts relate to one another. Here, all varieties are over a field of characteristic $0$ (the necessary condition for the below universality results).

The universality of $\mathbb{OB}$ among similarly behaved oriented bivariant theories for quasi-projective varieties was proved in \cite{sy} and \cite{y}. If $A$ is any such theory, we get a map of oriented bivariant theories $\mathbb{OB} \rightarrow A$, which in particular gives a group homomorphism $\mathbb{OB}(X \xrightarrow{f} Y) \rightarrow A(X \xrightarrow{f} Y)$ for each morphism $(X \xrightarrow{f} Y)$ in $\mathcal{V}$. We have then, for each such morphism, four group homomorphisms $\mathbb{OB}(X \xrightarrow{f} Y) \rightarrow \mathbb{OB}1(X \xrightarrow{f} Y)$, $\mathbb{OB}(X \xrightarrow{f} Y) \rightarrow \mathbb{OB}2(X \xrightarrow{f} Y)$, $\mathbb{OB}(X \xrightarrow{f} Y) \rightarrow \mathbb{OB}3(X \xrightarrow{f} Y)$ and $\mathbb{OB}(X \xrightarrow{f} Y) \rightarrow \mathbb{OB}4(X \xrightarrow{f} Y)$. This again suggests that the four previous bivariant theories are given as quotients of $\mathbb{OB}$, which is constructed from ${\mathbb M}(\mathcal{V}/X \xrightarrow{f} Y)$. \\

The constructions in the previous sections gave the following covariant and contravariant theories. Recall that each construction was made in order to generalize one of the four instances of algebraic cobordism; this is reflected in the group isomorphisms below. \\

$\mathbb{OB}1_*(X) := \mathbb{OB}1(X \xrightarrow{c} pt) \cong \Omega_*(X) \cong \omega_*(X)$ for $X$ quasi-projective. \\

$\mathbb{OB}1^*(X) := \mathbb{OB}1(X \xrightarrow{id} X)$  is a contravariant theory for $X$ smooth.\\

$\mathbb{OB}2_*(X) := \mathbb{OB}2(X \xrightarrow{c} pt)$ is a covariant theory for $X$ quasi-projective.\\

$\mathbb{OB}2^*(X) := \mathbb{OB}2(X \xrightarrow{id} X) \cong \Omega^*(X) \cong \omega^*(X)$ for $X$ smooth.\\

$\mathbb{OB}3_*(X) := \mathbb{OB}3(X \xrightarrow{c} pt) \cong \Omega_*(X) \cong \omega_*(X)$ for $X$ quasi-projective.\\

$\mathbb{OB}3^*(X) := \mathbb{OB}3(X \xrightarrow{id} X)$ is a contravariant theory for $X$ smooth. \\

$\mathbb{OB}4_*(X) := \mathbb{OB}4(X \xrightarrow{c} pt)$ is a covariant theory for $X$ quasi-projective.\\

$\mathbb{OB}4^*(X) := \mathbb{OB}4(X \xrightarrow{id} X) \cong \Omega^*(X) \cong \omega^*(X)$ for $X$ smooth.\\

Algebraic bordism $\Omega_*$ is universal among all oriented Borel-Moore functors $A_*$ with product and of geometric type for quasi-projective varieties (\cite{lm}, Thm. 2.4.13). This gives a natural transformation $\Omega_* \rightarrow A_*$ for each such $A_*$, and in particular a group homomorphism $\Omega_*(X) \rightarrow A_*(X)$ for each quasi-projective variety $X$.

This gives us group homomorphisms $\Omega_*(X) \rightarrow \mathbb{OB}2_*(X)$ and $\Omega_*(X) \rightarrow \mathbb{OB}4_*(X)$ for each quasi-projective variety $X$.

Finally, algebraic cobordism $\Omega^*(X)$ is universal among oriented cohomology theories $A^*$ for smooth varieties (\cite{lm}, Thm. 1.2.6), which gives a natural transformation $\Omega^* \rightarrow A^*$ for any such $A^*$, and in particular a ring homomorphism $\Omega^*(X) \rightarrow A^*(X)$ for each smooth variety $X$.

In our particular cases, we get ring homomorphisms $\Omega^*(X) \rightarrow \mathbb{OB}1^*(X)$ and $\Omega^*(X) \rightarrow \mathbb{OB}3^*(X)$ whenever $X$ is a smooth variety.

\end{document}